\documentclass{article}

\usepackage[german,english]{babel}
\usepackage{amsmath,amssymb,latexsym}
\usepackage{booktabs,tabularx}
\usepackage{graphicx,psfrag}
\usepackage{url}

\newtheorem{deff}{Definition}

\newtheorem{thm}[deff]{Theorem}
\newtheorem{cor}[deff]{Corollary}

\graphicspath{{.}}

\makeatletter

\title{Small Examples of Non-Constructible Simplicial Balls and Spheres}

\author{\Large Frank H.~Lutz}

\date{}

\begin{document}

\selectlanguage{english}

\maketitle

\medskip

\begin{abstract}
We construct non-constructible simplicial $d$-spheres with $d+10$ vertices
and non-constructible, non-realizable simplicial $d$-balls with $d+9$ vertices
for $d\geq 3$.
\end{abstract}

\section{Introduction}
The concepts of \emph{vertex-decomposability}, \emph{shellability},
and \emph{constructibility} describe three particular ways to assemble
a simplicial complex from the collection of its facets (cf.\ Bj\"orner \cite{Bjoerner1995}).
The following implications are strict for (pure) simplicial complexes:
\begin{center}
vertex decomposable\, $\Longrightarrow$\, shellable\, $\Longrightarrow$\, constructible.
\end{center}

Shellability has its origin in Schl\"afli's computation from 1852 \cite{Schlaefli1901} 
of the \mbox{Euler} characteristics of convex polytopes, where he
based his calculation on the assumption that the boundary complexes 
of polytopes are shellable. However, this property of polytopes was justified 
only much later in 1970 by Bruggesser and Mani \cite{BruggesserMani1971} and then played a crucial role
in McMullen's proof of the Upper Bound Theorem in the same year \cite{McMullen1970}. 
Besides in polyhedral theory, shellability has found fruitful applications in topology, 
combinatorics, and computational geometry; see the surveys \cite{Bing1964}, \cite{Bjoerner1995}, 
\cite{DanarajKlee1978}, \cite[Ch.~8]{Ziegler1995}, \cite{Ziegler1998},
and the references contained therein.

The notion of constructibility was coined by Hochster in 1972~\cite{Hochster1972},
but implicitly was used long before in combinatorial topology. 
In particular, it follows from Newman's and Alexander's fundamental works
on the foundations of combinatorial and PL topology from 
1926 \cite{Newman1926c} and 1930 \cite{Alexander1930} (cf.\ also Bj\"orner \cite{Bjoerner1995})
that a constructible $d$-dimensional simplicial complex in which every $(d-1)$-face 
is contained in exactly two or at most two $d$-dimensional facets 
is a PL $d$-sphere or a PL $d$-ball, respectively.
For recent surveys on constructibility see \cite{Hachimori1999} and \cite{HachimoriZiegler2000}.

The strongest concept, vertex-decomposability, was introduced by Provan and Billera
in their proof from 1980 \cite{ProvanBillera1980} that vertex decomposable 
simplicial complexes satisfy the simplicial form of the famous Hirsch conjecture (cf.\ \cite[p.~168]{Dantzig1963})
of linear programming.

Although boundary spheres of simplicial polytopes are shellable, 
Lockeberg~\cite{Lockeberg1977} constructed a simplicial $4$-polytope 
with $12$ vertices that is not vertex-de\-com\-pos\-able; 
and there even are not vertex-decomposable simplicial $4$-polytopes with
$10$ vertices \cite{KleeKleinschmidt1987} and not vertex-decomposable, 
non-polytopal simplicial $3$-spheres with $9$ vertices \cite{BokowskiBremnerLutzMartin2003pre}.
For two-dimensional balls and spheres it was proved by Bing~\cite{Bing1964}
that they are shellable and by Provan and Billera \cite{ProvanBillera1980}
that they are vertex-decomposable. Klee and Kleinschmidt
\cite{KleeKleinschmidt1987} also showed that all simplicial $d$-balls 
and all simplicial $d$-spheres with up to $d+3$, respectively $d+4$
vertices, are vertex-decomposable. However, for $d\geq 3$ there are not
vertex-decomposable simplicial $d$-balls with $d+4$ vertices 
and $10$ facets as well as not vertex-decomposable simplicial $d$-spheres 
with $d+6$ vertices; see \cite{BokowskiBremnerLutzMartin2003pre}
and \cite{Lutz2003epre}.

The first known example of a non-shellable cellular $3$-ball 
is due to Furch and appeared in 1924 \cite{Furch1924}. A non-shellable 
simplicial $3$-ball with $30$ vertices and $72$ facets was provided 
by Newman in 1926 \cite{Newman1926d}. 
Newman's ball is \emph{strongly non-shellable}, 
i.e., it has no \emph{free} facet that can be removed from
the triangulation without loosing ballness. 
Much smaller strongly non-shellable simplicial $3$-balls were obtained 
by Gr\"unbaum (cf.\ \cite{DanarajKlee1978}) with $14$ vertices and $29$ facets
and by Ziegler \cite{Ziegler1998} with $10$ vertices and $21$ facets. 
Rudin's $3$-ball \cite{Rudin1958} with $14$ vertices 
and $41$ tetrahedra gives a strongly non-shellable rectilinear triangulation 
of a tetrahedron with all the vertices on the boundary; the vertices
even can be moved slightly to yield a straight triangulation 
of a convex $3$-polytope with $14$ vertices \cite{ConnellyHenderson1980}. 
Ziegler's ball is realizable as a straight, yet non-convex ball in $3$-space. 
Coordinates for a rectilinear realization of Gr\"unbaum's ball 
can be found in \cite{Hachimori1999}. Vertex-minimal non-shellable
$3$-balls with $9$ vertices are enumerated in \cite{BokowskiBremnerLutzMartin2003pre};
see \cite{Lutz2003fpre} for a geometric realization of one of
these balls with $18$ facets.

The existence of non-constructible $3$-balls was shown
by Lickorish \cite{Lickorish1971} in 1971, but it remained unclear 
whether there are non-shellable $3$-spheres.
Non-shellable cell partitions of $S^3$ were first constructed
by Vince~\cite{Vince1985} in 1985 and then by Armentrout~\cite{Armentrout1994}.
In 1991, Lickorish \cite{Lickorish1991} described non-shellable triangulated 
$3$-spheres that contain a knotted triangle made of the sum of 
(at least) three trefoil knots.
In fact, is suffices to use one single trefoil knot:
\begin{thm} {\rm (Hachimori and Ziegler \cite{HachimoriZiegler2000})}\label{thm:hz}
If a triangulated $3$-ball or $3$-sphere contains \textbf{any} 
knotted triangle, then it is non-constructible (and thus non-shellable).
Moreover, a $3$-ball with a knotted spanning arc consisting
of at most $2$ edges is non-constructible.
\end{thm}
A first explicit, but large, non-constructible triangulated $3$-sphere with $f$-vector
$f=(381,2309,3856,1928)$ based on Furch's $3$-ball with a knotted 
spanning arc consisting of one edge was constructed by Hachimori \cite{Hachimori_url}.
Suspensions of such spheres produce non-constructible simplicial 
PL $d$-spheres in dimensions $d\ge 3$. 
Examples of small non-PL (and hence non-constructible) $d$-spheres 
of dimensions $d\ge 5$ with $d+13$ vertices can be found in \cite{BjoernerLutz2000}; 
see also \cite{BjoernerLutz2003}. Their construction makes use of
the double suspension theorem of Edwards \cite{Edwards1975}
(respectively of its generalization by Cannon~\cite{Cannon1979})
that double suspensions of non-spherical homology $d$-spheres give non-PL $(d+2)$-spheres.

\section{The Examples}

In the following, we employ the theorem of Hachimori and Ziegler
to construct simplicial PL $d$-spheres in dimensions $d\geq 3$ 
with only $d+10$ vertices that are non-constructible. 
From the enumeration in \cite{BokowskiBremnerLutzMartin2003pre} 
it follows that all $3$-spheres with \mbox{$n\leq 10$} vertices are shellable.
Hence, the non-constructible $3$-sphere $S^3_{13,56}$ with $13$ vertices
that we are going to obtain is, if not vertex-minimal, then 
close to vertex-minimality.

\begin{figure}
\begin{center}
\small
\psfrag{1}{1}
\psfrag{2}{2}
\psfrag{3}{3}
\psfrag{4}{4}
\psfrag{5}{5}
\psfrag{6}{6}
\psfrag{7}{7}
\psfrag{8}{8}
\psfrag{9}{9}
\psfrag{10}{10}
\psfrag{11}{11}
\psfrag{12}{12}
\includegraphics[width=.7\linewidth]{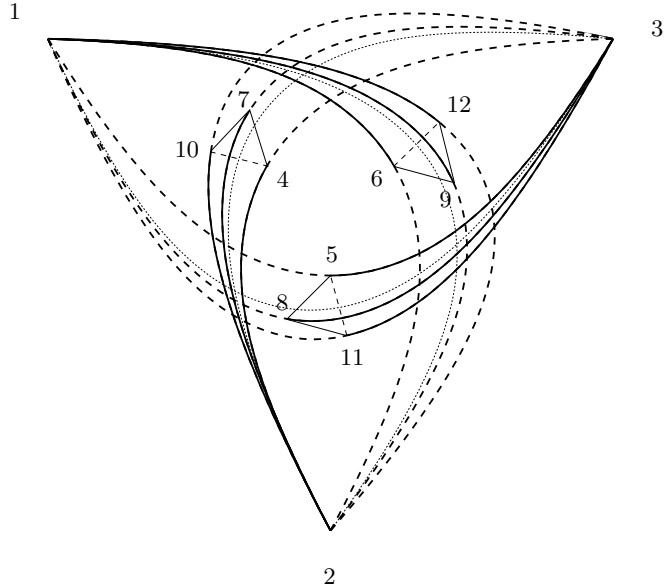}
\end{center}
\caption{The trefoil knot with three protected edges.}
\label{fig:trefoil}
\end{figure}

\begin{thm}
There is a non-constructible $3$-sphere $S^3_{13,56}$
with $13$ vertices and $56$ facets. Moreover, there are 
two strongly non-shellable, non-constructible $3$-balls $B^3_{12,37,a}$
and $B^3_{12,37,b}$ with $12$ vertices and $37$ facets that can not
be rectilinearly embedded into ${\mathbb R}^3$. 
\end{thm}
\emph{Proof.}
The examples are based on a trefoil knot consisting of three edges 
$12$, $13$, and $23$ (the dotted lines in Figure~\ref{fig:trefoil}) 
which we embed into ${\mathbb R}^3$. We shield off the edges 
by enclosing every edge with three tetrahedra,
as listed in the first column of Table~\ref{tbl:ball_3_16_46}.
We then close the holes of the knot by gluing in the following $16$ triangles:
\begin{center}\small
\begin{tabular}{l@{\hspace{10mm}}l@{\hspace{10mm}}l@{\hspace{10mm}}l}
$456$  &  $146$     &  $245$     &  $356$    \\
       &  $147$     &  $258$     &  $369$    \\
       &  $17\,10$  &  $28\,11$  &  $39\,12$ \\
       &  $15\,10$  &  $26\,11$  &  $34\,12$ \\
       &  $45\,10$  &  $56\,11$  &  $46\,12$.
\end{tabular}
\end{center}
\begin{table}[htp]
\small\centering
\defaultaddspace=0.3em
\caption{The ball $B^3_{16,46}$.}\label{tbl:ball_3_16_46}
\begin{center}\small
\begin{tabular*}{\linewidth}{@{\extracolsep{\fill}}llllll@{}}
\toprule
 \addlinespace
  $1269$     &  $146\,12$  &  $147\,13$     &  $258\,14$     &  $369\,15$     &  $456\,16$    \\
  $126\,12$  &  $245\,10$  &  $247\,13$     &  $358\,14$     &  $169\,15$     &  $146\,16$    \\
  $129\,12$  &  $356\,11$  &  $17\,10\,13$  &  $28\,11\,14$  &  $39\,12\,15$  &  $14\,13\,16$ \\
             &             &  $27\,10\,13$  &  $38\,11\,14$  &  $19\,12\,15$  &  $24\,13\,16$ \\
  $1358$     &             &  $15\,10\,13$  &  $26\,11\,14$  &  $34\,12\,15$  &  $245\,16$    \\
  $135\,11$  &             &  $25\,10\,13$  &  $36\,11\,14$  &  $14\,12\,15$  &  $25\,14\,16$ \\
  $138\,11$  &             &  $158\,13$     &  $269\,14$     &  $347\,15$     &  $35\,14\,16$ \\
             &             &  $258\,13$     &  $369\,14$     &  $147\,15$     &  $356\,16$    \\
  $2347$     &             &                &                &                &  $36\,15\,16$ \\
  $234\,10$  &             &                &                &                &  $16\,15\,16$ \\
  $237\,10$  &             &                &                &                & \\
 \addlinespace
\bottomrule
\end{tabular*}
\end{center}
\end{table}
The resulting simplicial complex $C$ is contractible.
By adding the $37$ tetrahedra in the columns 2--6 of
Table~\ref{tbl:ball_3_16_46} we thicken $C$
to a ball $B^3_{16,46}$ with $16$ vertices, $46$ facets,
and $f$-vector $f=(16,75,106,46)$.
\begin{figure}
\begin{center}
\small
\psfrag{1}{1}
\psfrag{2}{2}
\psfrag{3}{3}
\psfrag{4}{4}
\psfrag{5}{5}
\psfrag{6}{6}
\psfrag{7}{7}
\psfrag{8}{8}
\psfrag{9}{9}
\psfrag{10}{10}
\psfrag{11}{11}
\psfrag{12}{12}
\psfrag{13}{13}
\psfrag{14}{14}
\psfrag{15}{15}
\includegraphics[width=.7\linewidth]{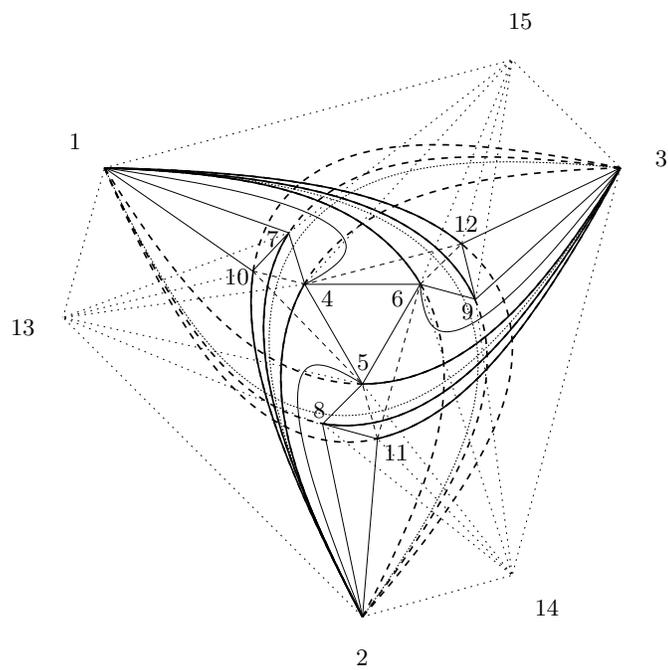}
\end{center}
\caption{The contractible complex $C$ with three cones.}
\label{fig:trefoil_16}
\end{figure}
Since $B^3_{16,46}$ contains a trefoil knot composed of three
edges, it follows from Theorem~\ref{thm:hz} of Hachimori and Ziegler
that $B^3_{16,46}$ is not constructible and thus not shellable.
In fact, $B^3_{16,46}$ is strongly non-shellable, as the removal 
of any of its facets destroys the ballness. Moreover, the presence 
of the $3$-edge knot prevents $B^3_{16,46}$ from having a straight 
embedding into ${\mathbb R}^3$. 

In Figure~\ref{fig:trefoil_16} we display the complex $C$.
We also indicate the cones with respect to the vertices $13$, $14$, and $15$
over eight of the triangles of $C$ each, as listed in columns 3--5 
of Table~\ref{tbl:ball_3_16_46}. The cone with respect to vertex $16$ 
is then placed ``above'' the drawing. 

The boundary of $B^3_{16,46}$ consists of $28$ triangles:
\begin{center}\small
\begin{tabular}{l@{\hspace{10mm}}l@{\hspace{10mm}}l@{\hspace{10mm}}l@{\hspace{10mm}}l}
$1\,13\,16$  &  $456$  &  $45\,10$  &  $56\,11$  &  $46\,12$ \\
$2\,13\,16$  &         &  $15\,10$  &  $26\,11$  &  $34\,12$ \\
$2\,14\,16$  &         &  $15\,11$  &  $26\,12$  &  $34\,10$ \\
$3\,14\,16$  &         &  $18\,11$  &  $29\,12$  &  $37\,10$ \\
$3\,15\,16$  &         &  $28\,11$  &  $39\,12$  &  $17\,10$ \\
$1\,15\,16$  &         &  $28\,13$  &  $39\,14$  &  $17\,15$ \\
             &         &  $18\,13$  &  $29\,14$  &  $37\,15$.
\end{tabular}
\end{center}
If we add to $B^3_{16,46}$ the cone over these $28$ triangles with
respect to a new vertex~$17$, then we get a $3$-sphere $S^3_{17,74}$ 
with $f=(17,91,148,74)$.
This $3$-sphere still contains the complex $C$ and with it
the trefoil knot composed of the three edges $12$, $13$, and~$23$.
Hence, $S^3_{17,74}$ is a not constructible, non-shellable sphere.
By construction, $B^3_{16,46}$ and $S^3_{17,74}$ have a
${\mathbb Z}_3$-symmetry.

\begin{figure}
\begin{center}
\small
\psfrag{1}{1}
\psfrag{2}{2}
\psfrag{3}{3}
\psfrag{4}{4}
\psfrag{5}{5}
\psfrag{6}{6}
\psfrag{7}{7}
\psfrag{8}{8}
\psfrag{9}{9}
\psfrag{10}{10}
\psfrag{11}{11}
\psfrag{12}{12}
\includegraphics[width=.7\linewidth]{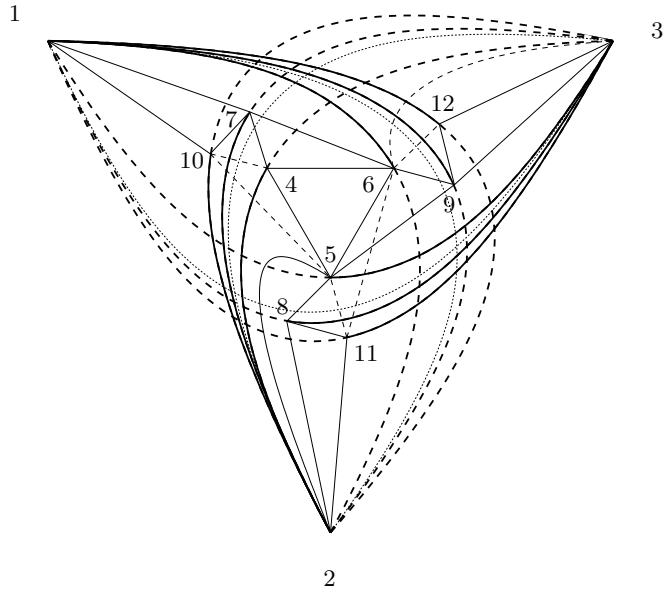}
\end{center}
\caption{The $3$-edge trefoil knot lying in the non-shellable sphere $S^3_{13,56}$.}
\label{fig:trefoil_12}
\end{figure}

Since all $3$-spheres with $n\leq 10$ vertices
are shellable \cite{BokowskiBremnerLutzMartin2003pre}, 
$17$ vertices is close to the minimal number of vertices 
that are needed for a non-shellable $3$-sphere.
In order to still improve on the number of vertices, we applied the
bistellar flip program BISTELLAR \cite{BIST} to $S^3_{17,74}$, 
under the additional restriction that the edges of the knot should not
be touched.
(The objective of BISTELLAR is to decrease the size of a triangulation
of a manifold by performing bistellar flips that locally modify
the triangulation without changing the topological type; 
see \cite{BjoernerLutz2000} for an explicit description.) 
As result, we obtained a simplicial $3$-sphere $S^3_{13,56}$ with $f=(13,69,112,56)$. 
The removal of the star of vertex $13$
\begin{center}\small
\begin{tabular}{l@{\hspace{10mm}}l@{\hspace{10mm}}l@{\hspace{10mm}}l@{\hspace{10mm}}l}
$179\,13$        &  $257\,13$     &  $358\,13$     &  $579\,13$        \\
$17\,11\,13$     &  $258\,13$     &  $359\,13$     &  $6\,10\,11\,13$  \\
$19\,10\,13$     &  $26\,11\,13$  &  $36\,10\,13$  \\
$1\,10\,11\,13$  &  $26\,12\,13$  &  $36\,12\,13$  \\
                 &  $27\,11\,13$  &  $38\,12\,13$  \\
                 &  $28\,12\,13$  &  $39\,10\,13$
\end{tabular}
\end{center}
from this complex yields a $12$-vertex $3$-ball $B^3_{12,38}$ with
$38$ facets, as listed in Table~\ref{tbl:ball_3_12_38}.
\begin{table}[htp]
\small\centering
\defaultaddspace=0.3em
\caption{The ball $B^3_{12,38}$.}\label{tbl:ball_3_12_38}
\begin{center}\small
\begin{tabular*}{\linewidth}{@{\extracolsep{\fill}}lllll@{}}
\toprule
 \addlinespace
$1269$     &  $158\,10$     &  \emph{2457} &  $3467$        &  $4567$       \\
$126\,12$  &  $15\,10\,11$  &  $245\,10$  &  \emph{346\,10} &  $456\,10$    \\ 
$129\,12$  &  $1679$        &  $258\,10$  &  $359\,11$     &  $5679$       \\
           &  $167\,12$     &  $269\,11$  &  $367\,12$     &  $569\,11$    \\  
$1358$     &  $178\,10$     &  $278\,10$  &  $37\,10\,12$  &  $56\,10\,11$ \\
$135\,11$  &  $178\,11$     &  $278\,11$  &  $389\,11$     \\
$138\,11$  &  $17\,10\,12$  &  $289\,11$  &  $389\,12$     \\
           &  $19\,10\,12$  &  $289\,12$  &  $39\,10\,12$  \\
$2347$     \\
$234\,10$  \\
$237\,10$  \\
 \addlinespace
\bottomrule
\end{tabular*}
\end{center}
\end{table}
This ball has two free facets, $2457$ and $346\,10$, so is not strongly
non-shellable. However, when we remove either of the two tetrahedra,
we get strongly non-shellable, non-constructible $3$-balls $B^3_{12,37,a}$
and $B^3_{12,37,b}$ with $37$ facets and $f=(12,58,84,37)$, respectively. 
These two balls are not isomorphic, although they have isomorphic boundaries. 
Both balls (and also the sphere $S^3_{13,56}$) still contain the original 
$3$-edge trefoil knot for which, this time, the triangles
\begin{center}\small
\begin{tabular}{l@{\hspace{10mm}}l@{\hspace{10mm}}l@{\hspace{10mm}}l}
$456$  &  $467$     &  $245$     &  $569$    \\
       &  $167$     &  $258$     &  $359$    \\
       &  $17\,10$  &  $28\,11$  &  $39\,12$ \\
       &  $15\,10$  &  $26\,11$  &  $36\,12$ \\
       &  $45\,10$  &  $56\,11$  &  $346$
\end{tabular}
\end{center}
are glued in to close the holes of the knot; see
Figure~\ref{fig:trefoil_12}. \hfill $\Box$

\begin{cor}
For $d\geq 3$ there are non-constructible $d$-spheres with $d+10$
vertices. Also there are non-constructible $d$-balls, $d\geq 3$, 
with $d+9$ vertices and $37$ facets that do not have a straight 
embedding into ${\mathbb R}^d$.
\end{cor}

\emph{Proof.}
The cone over a non-con\-structible, non-realizable $d$-ball
is a non-constructible, non-realizable $(d+1)$-ball with the 
same number of facets.
Similarly, the one-point suspension of a non-constructible $d$-sphere
is a non-constructible $(d+1)$-sphere; see \cite{JoswigLutz2003pre}.
\hfill $\Box$

\subsection*{Acknowledgment} 

The author is grateful to G\"unter M.~Ziegler
for helpful remarks.

\enlargethispage{3mm}

\bibliography{references_non_shellable}

\bigskip
\medskip


\noindent
\normalsize
Technische Universit\"at Berlin\\
Fakult\"at II - Mathematik und Naturwissenschaften\\
Institut f\"ur Mathematik, Sekr.\ MA 6-2\\
Stra\ss e des 17.\ Juni 136\\
10623 Berlin\\
Germany\\
{\tt lutz@math.tu-berlin.de}

\end{document}